\pdfoutput=1
\documentclass[sn-mathphys-num]{sn-jnl}

\usepackage{graphicx}%
\usepackage{multirow}%
\usepackage{amsmath,amssymb,amsfonts}%
\usepackage{amsthm}%
\usepackage{bm}
\usepackage{mathrsfs}%
\usepackage[title]{appendix}%
\usepackage{textcomp}%
\usepackage{manyfoot}%
\usepackage{booktabs}%
\usepackage{algorithm}%
\usepackage{algorithmicx}%
\usepackage{algpseudocode}%
\usepackage{listings}%

\usepackage{etoolbox}   
\makeatletter
\patchcmd{\@maketitle}
  {\artauthors}
  {\begin{center}\artauthors\end{center}}
  {}{}
\makeatother

\usepackage{mathtools}
\mathtoolsset{showonlyrefs}

\usepackage{algorithm}
\usepackage{algpseudocode}

\usepackage{color}

\usepackage{pgfplots}
\usetikzlibrary{intersections}
\usepgfplotslibrary{fillbetween}
\usepackage{filecontents}

\theoremstyle{thmstyleone}%
\theoremstyle{thmstyletwo}%
\newtheorem{remark}{Remark}%

\theoremstyle{thmstylethree}%

\newcommand{\bbR}{\mathbb{R}}

\newcommand{\bme}{\bm{e}}
\newcommand{\bmf}{\bm{f}}

\newcommand{\bmr}{\bm{r}}

\newcommand{\bmv}{\bm{v}}

\newcommand{\bmx}{\bm{x}}

\newcommand{\calO}{\mathcal{O}}

\newcommand{\rmd}{\mathrm{d}}
\newcommand{\rme}{\mathrm{e}}

\newcommand{\rmN}{\mathrm{N}}
\newcommand{\rmO}{\mathrm{O}}

\newcommand{\trans}{{\text{\tiny\sf T}}}

\newcommand{\Ga}{\mathrm{Ga}}
\newcommand{\GIG}{\mathrm{GIG}}

\newcommand{\iid}{\overset{\mathrm{i.i.d.}}{\sim}}

\begin{document}

\title[]{Quantifying uncertainty in the numerical integration of evolution equations based on Bayesian isotonic regression}

\author*[1]{\fnm{Yuto} \sur{Miyatake}}\email{yuto.miyatake.cmc@osaka-u.ac.jp}

\author[2]{\fnm{Kaoru} \sur{Irie}}

\author[3,4]{\fnm{Takeru} \sur{Matsuda}}


\affil[1]{\orgdiv{D3 Center}, \orgname{Osaka University}, 
\orgaddress{
\country{Japan}}}

\affil[2]{\orgdiv{Faculty of Economics}, \orgname{The University of Tokyo}, 
\orgaddress{
\country{Japan}}}

\affil[3]{\orgdiv{Graduate School of Information Science and Technology}, \orgname{The University of Tokyo}, 
\orgaddress{
\country{Japan}}}

\affil[4]{
\orgname{RIKEN Center for Brain Science}, 
\orgaddress{
\country{Japan}}}


\abstract{
This paper presents a new Bayesian framework for quantifying discretization errors in numerical solutions of ordinary differential equations.
By modelling the errors as random variables, we impose a monotonicity constraint on the variances, referred to as discretization error variances.
The key to our approach is the use of a shrinkage prior for the variances coupled with variable transformations. 
This methodology extends existing Bayesian isotonic regression techniques to tackle the challenge of estimating the variances of a normal distribution. 
An additional key feature is the use of a Gaussian mixture model for the $\log$-$\chi^2_1$ distribution, enabling the development of an efficient Gibbs sampling algorithm for the corresponding posterior.
}

\keywords{Discretization error quantification, Bayesian isotonic regression, shrinkage prior, Gaussian mixture model, Gibbs sampler}



\maketitle

\section{Introduction}\label{sec1}

Given the ordinary differential equation (ODE) model
\begin{equation}
    \label{eq:ode_model}
    \frac{\rmd}{\rmd t} \bmx (t) 
    = \bmf (\bmx (t) ), 
    \quad
    \bmx (t_0) = \bmx_0 \in V,
\end{equation}
where $V$ is an appropriate finite-dimensional space to which the solution $\bmx(t)$ belongs, and $\bmf: V \to V$ is assumed to be sufficiently regular,
we are concerned with quantifying the behaviour of discretization errors, i.e. the errors induced by the discretization.

Studying the behaviour of discretization errors theoretically has long been a major focus in numerical analysis~\cite{hw91,hn93}. 
However, the need for quantitatively evaluating error behaviour has surged recently.
For instance, when the ODE model involves unknown parameters to be estimated from noisy time-series data, approximating the solution to the ODE is often necessary. 
If this approximation lacks sufficient accuracy, the parameter estimation may become significantly biased (see, for example,~\cite{cg17}).
Note that 
achieving accurate numerical solutions within a reasonable computational timeframe, especially with limited computational resources, remains a challenge
\footnote{We give several situations.
(i) For partial differential equations with multi-dimensions in space, spatial and temporal mesh size constraints can limit accuracy.
(ii) Excessively accurate computation is reliable but not always recommended as it may consume large amounts of electricity.
(iii) Sensitivity to initial conditions affects not only numerical but also analytical solutions for chaotic systems. 
(iv) In long-term integrations, local discretization errors tend to accumulate, affecting the overall result. 
(v) Whilst most ODE solvers employ step-size control techniques with preset tolerances, these primarily manage local errors without ensuring global accuracy.}.
Thus, it becomes critical to evaluate and, where possible, reduce the bias. Developing methods to quantify the behaviour of numerical errors would be a significant step toward achieving this objective.

In pursuit of these objectives, various methodologies have been developed to quantify discretization errors and assess reliability using probabilistic or statistical frameworks. 
Examples of such approaches include ODE filters and smoothers~\cite{ks20,ss19,tk19,ts21}, as well as perturbative methods~\cite{ag20,cg17,ls19,ls22}, which fall under the category of probabilistic numerics. 
The key idea shared in those approaches is to consider a statistical model for the discretization error.
For an extensive review of probabilistic numerics, we refer the reader to~\cite{ho22}.

In addition, the authors of this paper have developed slightly different methods~\cite{mmm24,mm21,mm23arxiv}
based on `isotonic regression.'
These methods consider a statistical model for the discretization error and estimate key parameters, referred to as discretization error variances, using both time-series noisy observations and numerical approximations. 
The key in these methods is incorporating the known structure of discretization errors into inference rather than simply modelling the errors as independent random variables.
Since discretization errors tend to accumulate as time integration proceeds in many numerical problems,
these methods impose a monotonicity constraint on the discretization error variances. 
Our previous approaches focus on using maximum likelihood estimation for these variances, and the methods can be incorporated into inverse problems, where the parameters of the ODE model are estimated via the maximum likelihood method.

In this paper, we aim to develop a Bayesian framework for discretization error quantification, hoping that this approach can be integrated into inverse problems solved using Bayesian methods like \cite{oates2019bayesian}. 
We need to impose prior information for the discretization error variances so that the monotonicity assumption is met.
However, determining a suitable prior is challenging, as it cannot be uniquely specified.

In this study, we employ a shrinkage prior with certain variable transformations, building on a recent work on Bayesian isotonic regression~\cite{oh24}.
Unlike most works on Bayesian isotonic regression, which typically focus on estimating normal means~\cite{cd07,cg11,nd04,oh24,ss09}, our approach targets the estimation of the variances of a normal distribution.
We address this gap by approximating the $\log$-$\chi^2_1$ distribution with a Gaussian mixture model.
By incorporating this approximation, the conditional posteriors for all random variables (parameters) can be sampled from well-known distributions, allowing us to construct a highly efficient Gibbs sampling algorithm for the posterior distribution.

The paper is organized as follows.
Section~\ref{subsec:problem} provides a description of the problem.
Our new method for quantifying the discretization error is presented in Section~\ref{sec2}.
Then, the method is tested numerically in Section~\ref{sec3}.
Finally, concluding remarks are given in Section~\ref{sec4}.

\subsection{Description of the problem}
\label{subsec:problem}

To focus on the inferential problem of monotonically increasing error variances, we restrict our attention to the following setting throughout this paper. 
Assume that noisy time-series observations are obtained at discrete time points 
$t=t_1,t_2,\dots,t_n$ ($t_0\leq t_1 < t_2 < \dots < t_n$).
The observation operator is denoted by $\calO:V\to \bbR^p$, and the observation noise is assumed to be a $p$-dimensional Gaussian vector with mean zero and covariance matrix $\Gamma \in \bbR^{p\times p}$.
The observation at $t=t_i$ is denoted by $\bmv_i = \calO(\bmx(t_i)) + \bme_i$, where
$\bme_i \sim \rmN_p (\bm{0},\Gamma)$.
We assume that $\Gamma$ is given and diagonal.

We solve \eqref{eq:ode_model} numerically using some numerical integrators such as the Runge--Kutta method and the numerical approximation at $t=t_i$ is denoted by $\bmx_i \approx \bmx(t_i)$.
We aim to quantify the error $\calO(\bmx_i) - \calO(\bmx(t_i))$ for $i=1,\dots,n$, using the available numerical solutions and observations.
Note that our interest is not in estimating a \emph{strict} upper bound of the discretization error.

\vspace{11pt}
\begin{remark}
The observation time points $t = t_1, t_2, \dots, t_n$ are not necessarily equidistant. 
We do not necessarily assume that the numerical method employs a step size of $t_i - t_{i-1}$. 
Instead, a smaller step size is typically used.
If some of the $t_i$ values do not coincide with the computational grid points, we generally obtain a numerical approximation at $t = t_i$ by interpolating between the neighbouring numerical results.
\end{remark}

\section{New discretization error quantification}\label{sec2}

\subsection{Discretization error variances}

Building on our previous work~\cite{mmm24,mm21,mm23arxiv}, we model the discretization errors for the observed variables, $\calO(\bmx(t_i)) - \calO(\bmx_i)$ for $i=1,\dots ,n$, as independent Gaussian random variables. 
Equivalently, the numerical solution $\calO(\bmx_i)$ at time $t = t_i$ is modelled as,
\begin{align}
    \calO(\bmx_i)  \sim \rmN_p (\calO(\bmx(t_i)),V_i),
\end{align}
where $V_i$ is a positive semi-definite $p\times p$ matrix.
The residual, which represents the difference between the observed data and the numerical approximation, is given by
\begin{equation}\label{lik1}
    \bmr_i = \bmv_i - \calO(\bmx_i) = (\bmv_i - \calO(\bmx(t_i))) + (\calO(\bmx(t_i))-\calO(\bmx_i)) \sim \rmN_p (\bm{0},\Gamma + V_i),
\end{equation}
where we assumed that $\bmv_i - \calO(\bmx(t_i))\sim \rmN_p (\bm{0},\Gamma)$ and $(\calO(\bmx(t_i))-\calO(\bmx_i))\sim \rmN_p (\bm{0},V_i)$ are independent.
Note that this residual is computed by observation $\bmv_i$ and numerical solution $\calO(\bmx_i)$, hence observable. 
Letting $\Sigma_i = \Gamma + V_i$ and $\Sigma = (\Sigma_1,\dots,\Sigma_N)$,
the likelihood function is expressed as
\begin{equation*}
    p(\bmr;\Sigma) = 
    \sum_{i=1}^N \frac{1}{\sqrt{(2\pi)^p \det \Sigma_i}} \exp \Big( -\frac{1}{2} \bmr_i^\trans \Sigma_i^{-1} \bmr_i \Big).
\end{equation*}
This likelihood does not involve $\calO(\bmx(t_i))$, the true value of the (mapped) ODE solution. We do not estimate $\calO(\bmx(t_i))$; the parameter of our interest is $\Sigma _i$. This likelihood can also be obtained as the marginal likelihood when we set an improper prior $ p( \calO(\bmx(t_i)) ) = 1$ and marginalize it out in the likelihood computed by $(\bmv_i,\calO(\bmx_i))$. 

With a prior for $\Sigma$, denoted as $\pi(\Sigma)$, we compute the posterior distribution
\begin{equation*}
\pi(\Sigma \mid \bmr) \propto \pi(\Sigma) p(\bmr ;\Sigma).
\end{equation*}
In the following subsections, we discuss the choice of $\pi (\Sigma)$ and develop a sampling algorithm for the posterior computation.

For the remainder of this paper, we assume that each $\Sigma_i$ is a diagonal matrix. 
Note that, in the context of maximum likelihood estimation for $\Sigma$, this diagonality assumption has been relaxed~\cite{mmm24}, and we believe similar approaches should be incorporated into the Bayesian framework. 
However, we leave the exploration of a more general positive semi-definite $\Sigma_i$ for future work.

\vspace{11pt}
\begin{remark}
    The suitability of assuming Gaussian distributions (with mean zero) for the discretization error depends on both the specific equation and the chosen numerical integrator, and this assumption may not necessarily capture the full dynamics of discretization error propagation. 
    We employ this assumption primarily to simplify the analysis and facilitate the implementation of a sampling algorithm.
    However, it is important to note that the statistical model is not inherently limited to Gaussian distributions, and alternative distributions can be considered in principle.
\end{remark}

\subsection{Prior distribution $\pi(\Sigma)$}

The prior distribution $\pi(\Sigma)$ is a crucial factor in our model.
It should reflect the actual behaviour of the discretization error, which 
typically accumulates as the time integration proceeds.
To integrate the knowledge of this behaviour into inference, we consider a prior distribution that ensures the corresponding posterior satisfies the condition $\Gamma \preceq \Sigma_1 \preceq \Sigma_2 \preceq \dots \preceq \Sigma_n$. 

Accordingly, the probability density function of the prior must have non-negative values only in the region where $\Gamma \preceq \Sigma_1 \preceq \Sigma_2 \preceq \dots \preceq \Sigma_n$ and zero elsewhere.
Since we have assumed $\Sigma_i$ is diagonal, the condition above reduces to the monotonicity of each diagonal entries. Thus, without loss of generality, we focus on the univariate case where $\calO:V\to \bbR$ (i.e. $p=1$), and 
express $\Gamma = \gamma^2$, $\Sigma_i = \sigma_i^2$ and $\Sigma= (\sigma_1^2,\dots,\sigma_n^2)$.

Bayesian approaches incorporating the monotonicity assumption have been extensively studied, typically within the framework of Bayesian isotonic regression~\cite{cd07,cg11,nd04,oh24,ss09}. While most studies in the literature focus on enforcing monotonicity for normal means, our interest lies in the monotonicity of the variances. 
Despite this distinction, many of the priors proposed in the literature could, in principle, be adapted to our case. 
However, the direct application of these priors introduces challenges in sampling from the posterior distribution, which we will discuss later.

In the authors' previous work~\cite{mmm24,mm21}, isotonic regression was used for discretization error quantification, where the estimated variances frequently exhibited step function-like behaviour. 
This observation has motivated us to explore the use of a half-shrinkage prior, specifically the horseshoe prior, for the discretization error variances after a certain change of variables. 
Our prior is developed based on the recent work by \cite{oh24}.

First, we transform $\sigma_i^2$ as follows:
\begin{equation}
    \begin{split}
        & \eta_1 = \log \sigma_1^2, \\
    & \eta_j = \log \sigma_j^2 - \log \sigma_{j-1}^2, \quad j=2,\dots,n.
    \end{split}
    \label{eq:eta_def}
\end{equation}
Under this transformation, the constraint
\begin{equation}
    0<\gamma^2 \leq \sigma_1^2  \leq \sigma_2^2 \leq \dots \leq \sigma_n^2
\end{equation}
becomes
\begin{equation}
    \log \gamma^2 \leq \eta_1, \quad 0\leq \eta_j, \quad j = 2,\dots,n.  
\end{equation}
Below, we use the notation $\rmN_{\geq \alpha}(\mu,\sigma^2)$ 
to describe the truncated normal distribution with mean $\mu$ and variance $\sigma^2$ over the domain $[\alpha, \infty)$. 

\vspace{11pt}
 \begin{remark}
    The transformation \eqref{eq:eta_def}, together with the condition $\eta_j \ge 0$, ensures the required constraint, $\sigma _{j-1}^2 \le \sigma_j^2$. This is equivalent to bounding the ratio as $\sigma _j^2 / \sigma_{j-1}^2 \ge 1$, which is typical practice in statistical modelling. The utility of taking the log-transformation is to make $\eta_j$ be real-valued and allow for the use of (truncated) Gaussian models for $\eta_j$. 
    
    Instead of this transformation, one could also consider an alternative approach:
     \begin{equation}
     \begin{split}
         & \eta_1 =  \sigma_1^2, \\
     & \eta_j = \sigma_j^2 - \sigma_{j-1}^2, \quad j=2,\dots,n,
     \end{split}
     \end{equation}
     and enforce the monotonicity on these $\eta_1,\dots,\eta_n$. 
     However, the likelihood function becomes a complicated function of $\eta_j$, for which the sampling of $\eta_j$ becomes more challenging. 
\end{remark}

We define the shrinkage prior for $\eta_j$ by using the scale mixture of truncated-normals.  
Introducing several additional random variables, the prior for $\eta_j$ is conditionally Gaussian and expressed as follows.
\begin{itemize}
    \item Prior for $\eta_1$:
    \begin{equation*}
    \eta_1 \mid \tau_1 \sim \rmN_{\geq \log \gamma^2} (\log \gamma^2,\tau_1),
    \quad 
    \tau_1 \mid \nu_1 \sim \Ga (1,\nu_1),
    \quad 
    \nu_1 \sim \Ga (\tfrac12 , 1),
    \end{equation*}
    \item Priors for $\eta_j$ ($j=2,\dots,n$):
    \begin{align*}
    & \eta_j\mid \tau_j,\lambda \sim \rmN_{\geq 0}(0,\lambda \tau_j),\\
    & \tau_j \mid \nu_j \sim \Ga (\tfrac12,\nu_j), \quad 
    \nu_j \sim \Ga (\tfrac12 , 1), \\
    & \lambda \mid \xi \sim \Ga (\tfrac12,\xi), \quad \xi \sim \Ga (\tfrac12, 1).
    \end{align*}
\end{itemize}
Here, the density of Gamma distribution $\Ga (\alpha, \beta)$ with shape parameter $\alpha>0$ and rate parameter $\beta>0$ is given by
\begin{equation*}
    \Ga(x|\alpha ,\beta ) = \frac{\beta^\alpha}{\Gamma (\alpha)} x^{\alpha-1} \rme^{-\beta x}.
\end{equation*}

The priors for $\eta_j$ ($j=2,\dots,n$) play an essential role in our model.
We refer to $\tau_j$ as the local parameters, $\lambda$ as the global parameters, and the remaining variables as the latent parameters. 
The global parameter $\lambda$ shrinks all $\eta_j$'s towards zero uniformly, while the local parameter $\tau_j$ provides a customized shrinkage effect for each $j$.

With the latent variables being marginalized out, $\sqrt{\tau_j}$ and $\sqrt{\lambda}$ follow the half-Cauchy distribution \citep{car09,car10}. 
The reason for the explicit use of the latent parameters is to simplify the construction of a sampling algorithm from the posterior distribution, as we will explain later.

The latent parameters also provide intuitive insight into the roles of $\tau_j$ and $\lambda$.
For example, since $\tau_j \mid \nu_j \sim \Ga (\tfrac12,\nu_j)$, the density of $\tau_j$ diverges as $\tau_j\to 0$, meaning very small values of $\tau_j$ are taken with very high probability.
This suggests a strong shrinkage effect, where $\sigma_j^2$ and  $\sigma_{j-1}^2$ can be very close.
On the other hand, since $\nu_j \sim \Ga (\tfrac12 , 1)$, 
the density of $\tau_j$ does not decay rapidly as $\tau_j \to \infty$ (this can also be understood by noticing that $\sqrt{\tau_j}$ follows the half-Cauchy distribution).
As a result, relatively large values of $\tau_j$ can also be sampled occasionally, allowing the model to accommodate sharp increases in $\sigma_j^2$.

\subsection{Sampling from the posterior distribution }

To simplify the notation, we write the collection of parameters $\nu_j$'s, $\tau_j$'s, $\xi$ and $\lambda$ by $\Upsilon$. Also, we write the collection of the other parameters, $\sigma_j^2$'s, by $\Sigma$, but note that these parameters can be written as the function of $\eta_j$'s. 
Then, the posterior distribution of interest is $\pi(\Sigma, \Upsilon \mid r) = \pi(\Sigma, \Upsilon) p(r;\Sigma) $, 
where $r = (r_1,\dots,r_n)$ represents the residuals, with $r_i = v_i - \calO(\bmx_i)$.
The posterior density is explicitly written down as
\begin{align*}
    & \pi(\Sigma,\Upsilon \mid \bmr) \\
    & \quad \propto
    \prod_{i=1}^n \frac{1}{\sqrt{2\pi \sigma_i}} \exp \left( -\frac{1}{2\sigma_i^2} r_i^2 \right)
\cdot
    \tau_1^{-1/2} \exp \Big( -\frac{(\eta_1-\log \gamma^2)^2}{2\tau_1}\Big) \cdot
    \nu_1 \rme^{-\nu_1 \tau_1} \cdot \nu_1^{-1/2}\rme^{-\nu_1} \\
    & \quad \quad \phantom{\propto} \cdot
    \prod_{j=2}^n \lambda^{-1/2}
    \tau_j^{-1/2} \exp \Big( -\frac{\eta_j^2}{2\lambda \tau_j}\Big) \cdot
    \nu_j^{1/2} \tau_j^{-1/2}\rme^{-\nu_j \tau_j} \cdot \nu_j^{-1/2} \rme^{-\nu_j} \\
    & \quad \quad \phantom{\propto} \cdot
    \xi^{1/2} \lambda^{-1/2} \rme^{-\xi\lambda} \cdot \xi^{-1/2} \rme^{-\xi}.
\end{align*}
To evaluate the posterior distribution above, we employ Gibbs sampling, which requires the sampling of each random variable, conditional on the others, from its full conditional posterior distribution. 
The samples from the full posterior distribution are obtained by iteratively sampling from the full conditional distributions.

\subsubsection{Conditional posteriors for the local, global and latent variables}
First, we consider the full conditional posterior distributions for the local, global and latent parameters, or $\Upsilon$, 
which are shown to be the
gamma and generalized inverse Gaussian distributions as follows.
\begin{align*}
    & \nu_1 \mid \tau_1 \sim \Ga (\tfrac32, 1 + \tau_1), \\
    & \tau_1 \mid \eta_1, \nu_1 \sim \GIG (2\nu_1 ,(\eta_1-\log\gamma^2)^2 , \tfrac12), \\
    & \nu_j \mid \tau_j \sim \Ga (1,1+\tau_j), \\
    & \tau_j \mid \eta_j, \nu_j  \sim \GIG \bigg(2\nu_j, \frac{\eta_j^2}{\lambda}, 0\bigg), \\
    & \xi \mid \lambda \sim \Ga (1,1+\lambda), \\
    & \lambda \mid \eta , \xi \sim \GIG \bigg(2\xi, \sum_{j=2}^n\frac{\eta_j^2}{\tau_j},\frac{-n+2}{2}\bigg).
\end{align*}
Here, the density function of the generalized inverse Gaussian (GIG) distribution, $\GIG(a,b,p)$, with parameters $a>0$, $b>0$ and $p\in\bbR$, is given by
\begin{equation*}
    \GIG(x|a,b,p) = \frac{(a/b)^{p/2}}{2K_p(\sqrt{ab})}x^{p-1} \rme^{-(ax+b/x)/2}, \quad x>0,
\end{equation*}
where $K_p(\cdot)$ is the modified Bessel function of the second kind.

\vspace{11pt}
\begin{remark}
    The conventional sampling algorithm~\cite{da89,le89} for the GIG distribution is an acceptance-rejection method based on the Ratio-of-Uniforms method.
    However, we should treat the following two situations carefully.

    First, when $b$ approaches $0$, the GIG distribution resembles the Gamma distribution $\Ga(a,p)$. 
    This approximation is implemented in the R library when $b$ is close to the machine epsilon, and importantly, it can be justified within the framework of the independent Metropolis--Hastings method.\footnote{Suppose that, at one step of Gibbs sampler, we are required to sample from $\GIG(a,b,p)$ with extremely small $b$. Let $x^{\rm old}$ be the sample generated at the previous iteration. In the independent Metropolis-Hastings algorithm, we sample $x^{\rm new} \sim \Ga(p,a)$, then accept $x^{\rm new}$ as the sample of this iteration with probability $\alpha (x^{\rm new},x^{\rm old})$, otherwise reject $x^{\rm new}$ and take $x^{\rm old}$ as the sample of this iteration. The acceptance probability is computed as 
    \begin{equation*}
        \alpha (x^{\rm new} , x^{\rm old}  ) = \min \left\{ 1,\frac{ \GIG(x^{\rm new}|a,b,p) \Ga(x^{\rm old}|p,a) }{ \GIG(x^{\rm old}|a,b,p) \Ga(x^{\rm new}|p,a) } \right\}. 
    \end{equation*}
    }

    Second, when $p<1$, the performance of the conventional algorithm deteriorates if $|ab|$ is close to $0$.
    In such a case case, the algorithm proposed in~\cite{hl14} is recommended.
\end{remark}

\subsubsection{Conditional posteriors for the discretization error variances}

Next, we consider the sampling of $\eta_1,\dots,\eta_n$. 
Despite the log-transformation made in the definition of $\eta_j$'s, the full conditional posterior densities are still complex and difficult to directly sample from. To enable the sampling of $\eta _j$'s, we make a change of variables to $r_i$ so that $\eta_i$ is viewed as a location parameter, then approximate the non-Gaussian likelihood by the mixture of Gaussian distributions, so that the half-Gaussian prior for $\eta_j$ becomes conditionally conjugate. The former step of transformation originates from the study of stochastic volatility in finance, e.g.~\cite{nel88,har94}, and the latter step of approximation is known as mixture sampler in Bayesian econometrics \citep{kim98,omo07}.

We note that the residual $r_i$ can be expressed as
\begin{equation*}
    r_i = \varepsilon_i \exp (\eta_i/2), \quad \varepsilon_i \iid \rmN(0,1).
\end{equation*}
By defining $z_i = \log r_i^2$, we obtain
\begin{equation*}
    z_i = \eta_i + \tilde{\varepsilon}_i, \quad \tilde{\varepsilon}_i = \log \varepsilon_i^2.
\end{equation*}
The transformed error term, $\tilde{\varepsilon}_i$, follows the log-$\chi^2$ distribution. Its probability density function is given by
\begin{equation}
    \label{pdf:double-exponential}
    f(\tilde{\varepsilon}) = \frac{1}{\sqrt{2\pi}} \exp\Big( \frac{\tilde{\varepsilon}-\exp(\tilde{\varepsilon})}{2} \Big).
\end{equation}
Therefore, our likelihood can be rewritten as
\begin{equation}
    p(r;\Sigma) = \prod_{i=1}^n f(z_i - \eta_i).
\end{equation}

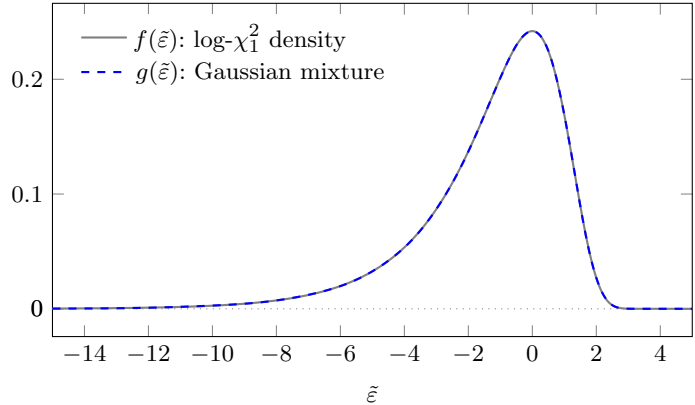
\begin{figure}
\centering
\begin{tikzpicture}
    \begin{axis}[
        width=10cm, 
        height=6cm,
        xlabel={$\tilde{\varepsilon}$},
        grid=none,            
        extra y ticks={0},           
        extra y tick style={
            grid=major,              
            grid style={thin, gray, dotted} 
        },
        ytick={0, 0.1, 0.2},     
        legend pos=north west, 
        xmin=-15, xmax=5,
        legend style={font=\small,align=left,draw=none},
        tick label style={font=\small},
        label style={font=\small},
    ]
        \addplot[mark=none, thick, gray] table [col sep=space] {gmm_data.dat};
        \addlegendentry{$f(\tilde{\varepsilon})$: log-$\chi^2_1$ density\ \ \ \ \ }
        \addplot[mark=none, thick, dashed, color=blue] table [col sep=space] {gmm1_data.dat};
        \addlegendentry{$g(\tilde{\varepsilon})$: Gaussian mixture}
    \end{axis}
\end{tikzpicture}
\caption{The probability distribution functions of the double exponential distribution $f(\tilde{\varepsilon})$ given in \eqref{pdf:double-exponential} and the Gaussian mixture model $g(\tilde{\varepsilon})$ give in \eqref{pdf:gmm}.}
\label{fig:gmm}
\end{figure}

The density function of $\tilde{\varepsilon}$ involves the double-exponential and complicates the likelihood, which hinders the direct sampling from the full conditionals of $\eta_j$'s. 
A practical remedy to this difficulty is to approximate $f$ using a Gaussian mixture model:
\begin{equation}
    \label{pdf:gmm}
    g(\tilde{\varepsilon}) = \sum_{k=1}^{10} w_k \rmN(\tilde{\varepsilon}\,;\, m_k, v_k^2).
\end{equation}
Mixture weight $w_k$, mixture component mean $m_k$ and variance $v_k^2$ are chosen so that the approximation error is minimized while matching the target and approximate first and second order moments. We used the values of $(w_k,m_k,v_k^2)$ computed in \cite{omo07}, as provided in Table~\ref{tab:omori07}. The probability density functions of $f(\tilde{\varepsilon})$ and $g(\tilde{\varepsilon})$ are displayed in Figure~\ref{fig:gmm}, in which it is clearly seen that the approximation of the log-$\chi_1^2$ distribution by the 10 Gaussian distributions is sufficiently accurate. 

\begin{table}[t]
    \centering
    \caption{Selection of $(w_j, m_j, v_j^2)$ given in~\cite{omo07}.}
    \label{tab:omori07}
    \begin{tabular}{c|rrr}
         $j$ & $w_j$ & $m_j$ & $v_j^2$ \\
         \hline
         $1$ & 0.00609 & 1.92677 & 0.11265 \\
         $2$ & 0.04775 & 1.34744 & 0.17788 \\
         $3$ & 0.13057 & 0.73504 & 0.26768 \\
         $4$ & 0.20674 & 0.02266 & 0.40611 \\
         $5$ & 0.22715 & $-0.85173$ & 0.62699 \\
         $6$ & 0.18842 & $-1.97278$ & 0.98583 \\
         $7$ & 0.12047 & $-3.46788$ & 1.57469\\
         $8$ & 0.05591 & $-5.55246$ & 2.54498 \\
         $9$ & 0.01575 & $-8.68384$ & 4.16591 \\
         $10$ & 0.00115 & $-14.65000$ & 7.33342\\
    \end{tabular}
\end{table}

\vspace{11pt}
\begin{remark}
    By taking the square of $\varepsilon_i$ to define $z_i$, we lose the information about the sign of residual $r_i$. This does not affect the inference, because the mean of $r_i$ is set to zero and is not estimated. When the mean of $r_i$ is modeled and estimated, or when $r_i$ is not conditionally independent of the other parameters, this loss of information must be addressed properly. In these situations, the approximating Gaussian mixture has been revisited and modified accordingly in \cite{hir24} and \cite{omo07}, respectively. 
\end{remark}

\vspace{11pt}
The mixture representation of $g(\cdot )$ is written in a hierarchical form with new latent parameters $s$ with $s_i \in \{ 1,2,\dots,10 \}$ for $i=1,\dots,n$. 
For each $i$, the hierarchical model defined by 
\begin{align}
    & \tilde{\varepsilon}_i | (s_i = k) \sim \rmN(\tilde{\varepsilon}\,;\, m_k, v_k^2), \\
    & p(s_i = k) = w_k, \qquad k \in \{ 1,\dots, 10\}, \label{eq:p_s}
\end{align}
has $g(\tilde{\varepsilon}_i )$ as the marginal distribution. 
We then propose the data augmentation, or sampling $(s_1,\dots, s_n)$ together with the other parameters. With $(s_1,\dots, s_n)$ being conditioned, the conditional posterior distributions of $\eta_1,\dots,\eta_n$ become tractable as follows.

\begin{itemize}
    \item Conditional posterior distribution for $s$:
    
    Since the latent parameters $s$ are independent of other parameters.
    $s_i$ can be sampled using the uniform distribution with the probability mass \eqref{eq:p_s}.
    
    \item Conditional posterior distribution for $\eta_1$, i.e. $\eta_1 \mid s, \tau_1$:

    Let 
    \begin{equation*}
    k_i := z_i - m_{s_i} - \eta_2 - \dots - \eta_{i}.
    \end{equation*}
    The conditional posterior for $\eta_1$ is found to be
    \begin{equation*}
        \eta_1 \mid s, \tau_1 \sim \rmN_{\geq \log\gamma^2} (\mu_1, w_1^2),
    \end{equation*}
    where 
    \begin{equation*}
        \mu_1 = \frac{\tau_1 \tilde{\mu}_1 + \tilde{\mu}_1^2 \log \gamma^2}{\tilde{w}_1^2 + \tau_1},
        \quad 
        w_1^2 = \frac{\tilde{w}_1^2 \tau_1}{\tilde{w}_1^2 + \tau_1}
    \end{equation*}
    with
    \begin{align*}
    & \tilde{\mu}_1 = \frac{\sum_{i=1}^n v_{s_1}^2\cdots v_{s_{i-1}}^2 v_{s_{i+1}}^2 \cdots v_{s_n}^2 k_i}{\sum_{i=1}^n v_{s_1}^2\cdots v_{s_{i-1}}^2 v_{s_{i+1}}^2 \cdots v_{s_n}^2}
    =
    \frac{\sum_{i=1}^n \frac{k_i}{v_{s_i}^2}}{\sum_{i=1}^n \frac{1}{v_{s_i}^2}},\\
    & \tilde{w}_1^2 = \frac{v_{s_1}^2 v_{s_2}^2\dots v_{s_n}^2}{\sum_{i=1}^n v_{s_1}^2\cdots v_{s_{i-1}}^2 v_{s_{i+1}}^2 \cdots v_{s_n}^2}
    =
    \frac{1}{\sum_{i=1}^n \frac{1}{v_{s_i}^2}}.
\end{align*}
    \item Conditional posterior distributions for $\eta_j$ ($j=2,\dots,n)$, i.e. $\eta_j \mid s, \lambda, \tau_j$:

    Let
    \begin{align*}
    k_{i,j} :&= z_i - m_{s_i} - \eta_1 - \dots - \eta_{j-1} - \eta_{j+1} - \dots - \eta_i \\
    &=k_{i,j-1} - \eta_{j-1} + \eta_j. 
    \end{align*}
    The conditional posteriors for $\eta_2,\dots,\eta_n$ are found to be
    \begin{equation*}
        \eta_j \mid s, \lambda, \tau_j \sim \rmN_{\geq 0}(\mu_j , w_j^2),
    \end{equation*}
    where
    \begin{align*}
    & \mu_j = \frac{\lambda \tau_j \tilde{\mu}_j}{\tilde{w}_j^2 + \lambda \tau_j}, \\
    & w_j^2 = \frac{\tilde{w}_j^2 \lambda \tau_1}{\tilde{w}_j^2 + \lambda \tau_j}
    \end{align*}
    with
    \begin{align*}
    & \tilde{\mu}_j = \frac{\sum_{i=j}^n v_{s_j}^2\cdots v_{s_{i-1}}^2 v_{s_{i+1}}^2 \cdots v_{s_n}^2 k_{i,j}}{\sum_{i=j}^n v_{s_j}^2\cdots v_{s_{i-1}}^2 v_{s_{i+1}}^2 \cdots v_{s_n}^2}
    =
    \frac{\sum_{i=j}^n \frac{k_{i,j}}{v_{s_i}^2}}{\sum_{i=j}^n \frac{1}{v_{s_i}^2}},\\
    & \tilde{w}_j^2 = \frac{v_{s_j}^2 v_{s_{j+1}}^2\cdots v_{s_n}^2}{\sum_{i=j}^n v_{s_j}^2\cdots v_{s_{i-1}}^2 v_{s_{i+1}}^2 \cdots v_{s_n}^2}
    =
    \frac{1}{\sum_{i=j}^n \frac{1}{v_{s_i}^2}}.
\end{align*}
\end{itemize}

Thus, all the full conditional distributions are shown to be well-known distributions, from which it is easy to simulate, enabling the use of Gibbs sampler.  
The overall Gibbs sampling algorithm is 
summarized in Algorithm~\ref{alg:gibbs_with_residual}.

\begin{algorithm}
\caption{Gibbs sampling algorithm for estimating the discretization error variances $\Sigma$}
\label{alg:gibbs_with_residual}
\begin{algorithmic}[1]
\State Input $z$ ($z_i = \log r_i^2$ for $i=1,\dots,n$)
\State Initialize $s^{(0)}, \eta^{(0)}, \nu^{(0)}, \tau^{(0)}, \xi^{(0)}, \lambda^{(0)}$
\For{$t = 1$ to $T$}
    \State Sample $s^{(t)} \sim P(s)$
    \State Sample $\eta_1^{(t)} \sim P (\eta_1 \mid s^{(t)}, \tau_1^{(t-1)})$
    \State Sample $\eta_j^{(t)} \sim P(\eta_j \mid s^{(t)}, \lambda ^{(t-1)}, \tau_j^{(t-1)})$ for $j=2,\dots,n$
    \State Sample $\nu_1^{(t)} \sim P(\nu_1 \mid \tau_1^{(t-1)})$
    \State Sample $\tau_1^{(t)} \sim P(\tau_1 \mid \eta_1^{(t)}, \nu_1^{(t)})$
    \State Sample $\nu_j^{(t)} \sim P(\nu_j \mid \tau_j^{(t-1)})$ for $j=2,\dots,n$
    \State Sample $\tau_1^{(t)} \sim P(\tau_j \mid \eta_j^{(t)}, \nu_j^{(t)})$ for $j=2,\dots,n$
    \State Sample $\xi^{(t)} \sim P(\xi\mid \lambda^{(t-1)})$
    \State Sample $\lambda^{(t)} \sim P(\lambda \mid \eta^{(t)},\xi^{(t)})$
    \State Compute $\{\sigma_i^2 \}^{(t)} = \exp(\eta_1^{(t)}+ \dots + \eta_i^{(t)})$
\EndFor
\State \Return $\Sigma^{(t)} = [\{\sigma_1^2\}^{(t)},\dots, \{\sigma_n^2\}^{(t)}]$ 
\end{algorithmic}
\end{algorithm}

\vspace{11pt}
\begin{remark}
    During the overall Gibbs sampling process,
    we must compute $k_{i,j}$ efficiently.
    To reduce the computational cost, 
    the denominator calculations for $\tilde{\mu}_j$ and $\tilde{w}_k^2$ should be implemented recursively.
    The computational complexity within a single loop is $\rmO (n^2)$.
    If the observed variable has the dimension $p$, and the sampling algorithm is applied across all observed variables, the total complexity per loop is $\rmO(pn^2)$.
    Note that it can be reduced to $\rmO(n^2)$ if implemented using a parallel architecture.

    While the complexity of numerically solving the ODE model is proportional to $n$, making the $\rmO(n^2)$ cost appear excessive by comparison, 
    it is important to note that the dimension of the ODE model is typically very large.
    In practical applications, this dimension can be in the millions or billions, especially when dealing with ODE models derived by discretizing partial differential equations.
    The computational cost per time step is then proportional to the model’s dimension, potentially scaling quadratically or cubically.
    Thus, when $n$ is much smaller than the ODE model's dimension, the $\rmO(n^2)$ cost of the sampling process might be acceptable.

    In our numerical experiments, we intentionally set the dimension of the ODE system to be relatively small to better observe the behaviour of the proposed approach. 
\end{remark}

\section{Numerical experiments} \label{sec3}

We test the proposed method on two systems: the FitzHugh--Nagumo (FN) model and the Kepler equation, both of which exhibit nearly periodic solutions.

For the FN model, we assume that only one of the dependent variables is observed. This scenario allows us to explore how the method handles partial observations and assesses the discretization error with a linear observation operator.

For the Kepler equation, we consider a scenario where the observation operator is a function of the dependent variables.
This setup introduces additional complexity, as the observations are influenced by multiple variables, providing a more rigorous test of the method's capability to quantify discretization errors under complex observation structures.

We perform the following numerical experiments using Julia v1.9.
The reference solutions, which are used instead of the exact solutions, are generated by the Tsitouras 5/4 Runge--Kutta method (\textsf{Tsit5()} with $\mathsf{abstol}=\mathsf{reltol} = 10^{-8}$ in the package DifferentialEquations.jl).
For comparison, we also show the results obtained by our previous approach~\cite{mm21} (the code is available\footnote{\url{https://github.com/yutomiyatake/IsoFuns.jl}}).
For the Gibbs sampling, we generated 2500 posterior samples after discarding 500 samples as a burn-in period. In predictive analysis, we increased the MCMC size to 29500 (with 500 burn-in) for the accurate estimation of posterior predictive tails. For the Kepler model, we generated 9500 samples after 500 burn-in. 

\subsection{The FitzHugh--Nagumo model}

\begin{figure}
    \centering
    \begin{tikzpicture}
    \begin{axis}[
        width=6cm, 
        xlabel={$t$},
        ylabel={$V$},
        ylabel style={yshift=-1em},
        legend pos=south west, 
        xmin=5, xmax=50,
        legend style=
        {fill opacity=0.8, font=\small,align=left},
        tick label style={font=\small},
        label style={font=\small},
    ]
        \addplot[mark=none, gray] table [col sep=space] {FN/refsol_v.dat};
        \addlegendentry{exact}

        \addplot[mark=*, only marks, mark size=1pt, gray] table [col sep=space] {FN/obs_v.dat};
        \addlegendentry{observation}

        \addplot[mark=none, thick, black] table [col sep=space] {FN/numsol_v.dat};
        \addlegendentry{numer. approx.}
    \end{axis}
\end{tikzpicture}
\begin{tikzpicture}
    \begin{axis}[
        width=6cm, 
        xlabel={$t$},
        ylabel={$R$},
        ylabel style={yshift=-1em},
        legend pos=south west, 
        xmin=5, xmax=50,
        legend style={fill opacity=0.8,font=\small,align=left},
        tick label style={font=\small},
        label style={font=\small},
    ]
        \addplot[mark=none, gray] table [col sep=space] {FN/refsol_r.dat};
        \addlegendentry{exact}

        \addplot[mark=*, only marks, mark size=1pt, gray] table [col sep=space] {FN/obs_r.dat};
        \addlegendentry{observation}

        \addplot[mark=none, thick, black] table [col sep=space] {FN/numsol_r.dat};
        \addlegendentry{numer. approx.}
    \end{axis}
\end{tikzpicture}

    \caption{Exact solution, observations and numerical approximations to the FN model.}
    \label{fig:FN-sol}
\end{figure}
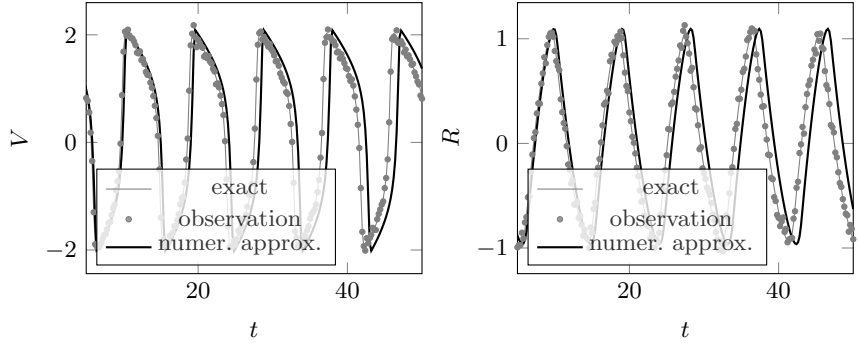

We consider the FitzHugh--Nagumo model
\renewcommand{\arraystretch}{1.1}
\begin{equation}
    \frac{\rmd}{\rmd t}
    \begin{bmatrix}
        V \\ R
    \end{bmatrix}
    =
    \begin{bmatrix}
        c \Big( V - \dfrac{V^3}{3} + R \bigg) \\
        -\dfrac{1}{c}(V-a+bR)
    \end{bmatrix},
    \quad
    \begin{bmatrix}
        V(0) \\ R(0) 
    \end{bmatrix}
    =
    \begin{bmatrix}
        1 \\ -1
    \end{bmatrix},
\end{equation}
where $(a,b,c) = (0.2,0.2,3)$.
We consider the two cases: $\calO([V,R]) = V$ or $\calO([V,R]) = R$.
For both cases, the observation noise variance is set to $\gamma^2 = 0.05^2$.
Each variable is observed from $t=5$ to $t=50$ with the interval $\Delta t = 0.2$.
As a numerical integrator, we employ the explicit Euler method with the step size $h=0.1$.
The exact solutions, observations and numerical approximations are plotted in Fig.~\ref{fig:FN-sol}.
It is observed that the propagation of the numerical solution is slower than the exact flow.

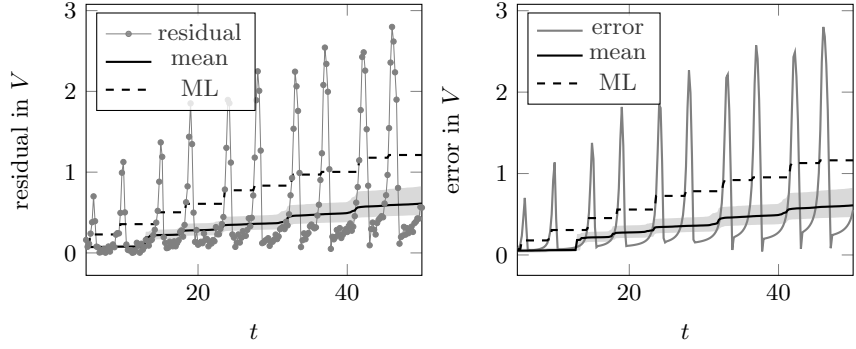
\begin{figure}
    \centering
    \begin{tikzpicture}
    \begin{axis}[
        width=6cm, 
        xlabel={$t$},
        ylabel={residual in $V$},
        ylabel style={yshift=-1em},
        legend pos=north west,
        xmin=5, xmax=50,
        legend style={fill opacity=0.8, font=\small,align=left},
        tick label style={font=\small},
        label style={font=\small},
    ]
        \addplot[mark=*,  mark size=1pt, line width=0.1pt, gray] table [col sep=space] {FN/res_v.dat};
        \addlegendentry{residual}

        \addplot[mark=none, thick, black] table [col sep=space] {FN/mean_v.dat};
        \addlegendentry{mean}

        \addplot[mark=none, thick, dashed] table [col sep=space] {FN/iso_v.dat};
        \addlegendentry{ML}

        \addplot[name path =upper, draw=none] table [col sep=space] {FN/upper_v.dat};
        \addplot[name path =lower,draw=none] table [col sep=space] {FN/lower_v.dat};
        \addplot[fill opacity=0.3, fill=gray] 
        fill between[of=upper and lower];
    \end{axis}
\end{tikzpicture}
\begin{tikzpicture}

    \begin{axis}[
        width=6cm, 
        xlabel={$t$},
        ylabel={error in $V$},
        ylabel style={yshift=-1em},
        legend pos=north west, 
        xmin=5, xmax=50,
        legend style={fill opacity=0.8, font=\small,align=left},
        tick label style={font=\small},
        label style={font=\small},
    ]
        \addplot[thick, gray] table [col sep=space] {FN/err_v.dat};
        \addlegendentry{error}

        \addplot[mark=none, thick, black] table [col sep=space] {FN/err_mean_v.dat};
        \addlegendentry{mean}

        \addplot[mark=none, thick, dashed] table [col sep=space] {FN/err_iso_v.dat};
        \addlegendentry{ML}

        \addplot[name path =upper, draw=none] table [col sep=space] {FN/err_upper_v.dat};
        \addplot[name path =lower,draw=none] table [col sep=space] {FN/err_lower_v.dat};
        \addplot[fill opacity=0.3, fill=gray] 
        fill between[of=upper and lower];
    \end{axis}
\end{tikzpicture}

    \caption{Quantification results for $V$ in the FN model.
    The left figure shows the absolute value of the residual $r_i$, the quantification results in terms of $\sigma_i$ with their mean and $95\%$ credible interval, along with the quantification based on the maximum likelihood approach~\cite{mm21}.
    The right figure shows essentially the same information but with the observation removed: the error is plotted instead of the residual, and $\sqrt{\sigma_i^2 - \gamma^2}$ is plotted for the quantification results.}
    \label{fig:FN-V-result}
\end{figure}

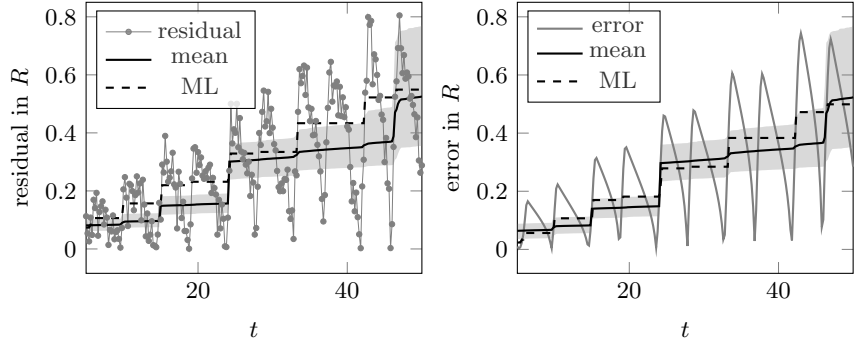
\begin{figure}
    \centering
    \begin{tikzpicture}
    \begin{axis}[
        width=6cm, 
        xlabel={$t$},
        ylabel={residual in $R$},
        ylabel style={yshift=-1em},
        legend pos=north west, 
        xmin=5, xmax=50,
        ymax=0.85,
        legend style={fill opacity=0.8, font=\small,align=left},
        tick label style={font=\small},
        label style={font=\small},
    ]
        \addplot[mark=*,  mark size=1pt, line width=0.1pt, gray] table [col sep=space] {FN/res_r.dat};
        \addlegendentry{residual}

        \addplot[mark=none, thick, black] table [col sep=space] {FN/mean_r.dat};
        \addlegendentry{mean}

        \addplot[mark=none, thick, dashed] table [col sep=space] {FN/iso_r.dat};
        \addlegendentry{ML}

        \addplot[name path =upper, draw=none] table [col sep=space] {FN/upper_r.dat};
        \addplot[name path =lower,draw=none] table [col sep=space] {FN/lower_r.dat};
        \addplot[fill opacity=0.3, fill=gray] 
        fill between[of=upper and lower];
    \end{axis}
\end{tikzpicture}
\begin{tikzpicture}

    \begin{axis}[
        width=6cm, 
        xlabel={$t$},
        ylabel={error in $R$},
        ylabel style={yshift=-1em},
        legend pos=north west,
        xmin=5, xmax=50,
        ymax=0.85,
        legend style={fill opacity=0.8, font=\small,align=left},
        tick label style={font=\small},
        label style={font=\small},
    ]
        \addplot[thick, gray] table [col sep=space] {FN/err_r.dat};
        \addlegendentry{error}

        \addplot[mark=none, thick, black] table [col sep=space] {FN/err_mean_r.dat};
        \addlegendentry{mean}

        \addplot[mark=none, thick, dashed] table [col sep=space] {FN/err_iso_r.dat};
        \addlegendentry{ML}

        \addplot[name path =upper, draw=none] table [col sep=space] {FN/err_upper_r.dat};
        \addplot[name path =lower,draw=none] table [col sep=space] {FN/err_lower_r.dat};
        \addplot[fill opacity=0.3, fill=gray] 
        fill between[of=upper and lower];
    \end{axis}
\end{tikzpicture}

    \caption{Quantification results for $R$ in the FN model.
    The left figure shows the residual $r_i$, the quantification results in terms of $\sigma_i$ with their mean and $95\%$ credible interval, along with the quantification based on the maximum likelihood approach~\cite{mm21}.
    The right figure shows essentially the same information but with the observation removed: the error is plotted instead of the residual, and $\sqrt{\sigma_i^2 - \gamma^2}$ is plotted for the quantification results.}
    \label{fig:FN-R-result}
\end{figure}

The quantification results for $V$ and $R$ are displayed in Fig.~\ref{fig:FN-V-result} and Fig.~\ref{fig:FN-R-result}, respectively.
For both cases, 
the absolute values of the residual and error are plotted.
The residual is defined as $r_i = v_i - x_i$, where $v_i$ is a noisy observation of $V$ or $R$, and $x_i$ is the numerical approximation, i.e. $x_i=V_i$ or $R_i$.
The error is given by $x(t_i) - x_i$, where $x(t_i) = V(t_i)$ or $R(t_i)$.
The left and right figures are essentially the same; however,
on the right figures, the observation is removed. 
Accordingly, 
in the left figures, the posterior mean and 95$\%$ credible interval are displayed based on 2,500 samples of $\sigma_i$.
For the right figures, similar results are plotted for
$\sqrt{\sigma_i^2 - \gamma^2}$.

It is observed that, for $V$, the proposed method yields smaller values than the maximum likelihood approach,
though the difference is at most about a factor of two.
Although these results may vary with changes in the observation time interval, this pattern is generally consistent.
In comparison, for $R$
the results of the proposed approach are closer to those of the maximum likelihood method.

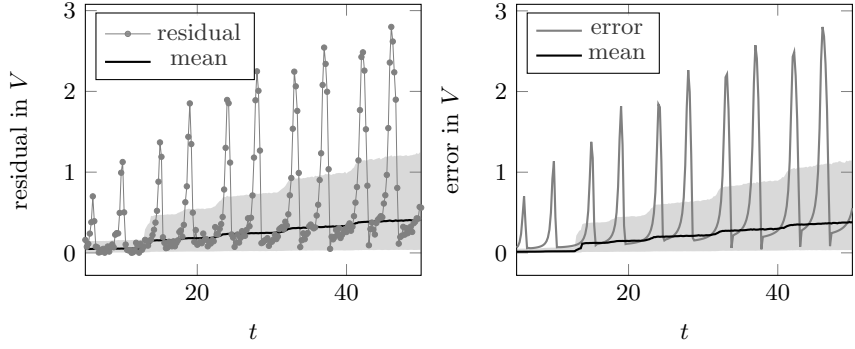
\begin{figure}
    \centering
    \begin{tikzpicture}
    \begin{axis}[
        width=6cm, 
        xlabel={$t$},
        ylabel={residual in $V$},
        ylabel style={yshift=-1em},
        legend pos=north west, 
        xmin=5, xmax=50,
        legend style={fill opacity=0.8, font=\small,align=left},
        tick label style={font=\small},
        label style={font=\small},
    ]
        \addplot[mark=*,  mark size=1pt, line width=0.1pt, gray] table [col sep=space] {FN/res_v.dat};
        \addlegendentry{residual}

        \addplot[mark=none, thick, black] table [col sep=space] {FN/res_sample_mean_v.dat};
        \addlegendentry{mean}

        \addplot[name path =upper, draw=none] table [col sep=space] {FN/res_sample_upper_v.dat};
        \addplot[name path =lower,draw=none] table [col sep=space] {FN/res_sample_lower_v.dat};
        \addplot[fill opacity=0.3, fill=gray] 
        fill between[of=upper and lower];
    \end{axis}
\end{tikzpicture}
\begin{tikzpicture}

    \begin{axis}[
        width=6cm, 
        xlabel={$t$},
        ylabel={error in $V$},
        ylabel style={yshift=-1em},
        legend pos=north west, 
        xmin=5, xmax=50,
        legend style={fill opacity=0.8, font=\small,align=left},
        tick label style={font=\small},
        label style={font=\small},
    ]
        \addplot[thick, gray] table [col sep=space] {FN/err_v.dat};
        \addlegendentry{error}

        \addplot[mark=none, thick, black] table [col sep=space] {FN/err_sample_mean_v.dat};
        \addlegendentry{mean}

        \addplot[name path =upper, draw=none] table [col sep=space] {FN/err_sample_upper_v.dat};
        \addplot[name path =lower,draw=none] table [col sep=space] {FN/err_sample_lower_v.dat};
        \addplot[fill opacity=0.3, fill=gray] 
        fill between[of=upper and lower];
    \end{axis}
\end{tikzpicture}

    \caption{
    These figures also present the residual and error plots for $V$. 
    In the left figure, 
    we sample $r_i$ from $\rmN(0,\sigma_i^2)$, where $\sigma_i^2$ is a sampled value using our Gibbs sampler, and 
    plot the $90\%$ credible interval for $|r_i|$.
    Similarly, in the right figure, we sample the error $\xi_i$ from $\rmN(0,\sigma_i^2 - \gamma^2)$ and plot the $90\%$ credible interval for $|\xi_i|$.
    For this test, we generated 29,500 posterior samples; however, even with a much smaller sample size, the results were almost identical (posterior means showed only slight oscillations).
    }
    \label{fig:FN-V-sample-result}
\end{figure}

\begin{figure}
    \centering
    \begin{tikzpicture}
    \begin{axis}[
        width=6cm, 
        xlabel={$t$},
        ylabel={residual in $R$},
        ylabel style={yshift=-1em},
        legend pos=north west,
        xmin=5, xmax=50,
        ymax=0.85,
        legend style={fill opacity=0.8, font=\small,align=left},
        tick label style={font=\small},
        label style={font=\small},
    ]
        \addplot[mark=*,  mark size=1pt, line width=0.1pt, gray] table [col sep=space] {FN/res_r.dat};
        \addlegendentry{residual}

        \addplot[mark=none, thick, black] table [col sep=space] {FN/res_sample_mean_r.dat};
        \addlegendentry{mean}

        \addplot[name path =upper, draw=none] table [col sep=space] {FN/res_sample_upper_r.dat};
        \addplot[name path =lower,draw=none] table [col sep=space] {FN/res_sample_lower_r.dat};
        \addplot[fill opacity=0.3, fill=gray] 
        fill between[of=upper and lower];
    \end{axis}
\end{tikzpicture}
\begin{tikzpicture}

    \begin{axis}[
        width=6cm, 
        xlabel={$t$},
        ylabel={error in $R$},
        ylabel style={yshift=-1em},
        legend pos=north west, 
        xmin=5, xmax=50,
        ymax=0.85,
        legend style={fill opacity=0.8, font=\small,align=left},
        tick label style={font=\small},
        label style={font=\small},
    ]
        \addplot[thick, gray] table [col sep=space] {FN/err_r.dat};
        \addlegendentry{error}

        \addplot[mark=none, thick, black] table [col sep=space] {FN/err_sample_mean_r.dat};
        \addlegendentry{mean}

        \addplot[name path =upper, draw=none] table [col sep=space] {FN/err_sample_upper_r.dat};
        \addplot[name path =lower,draw=none] table [col sep=space] {FN/err_sample_lower_r.dat};
        \addplot[fill opacity=0.3, fill=gray] 
        fill between[of=upper and lower];
    \end{axis}
\end{tikzpicture}

    \caption{
    These figures also present the residual and error plots for $R$. 
    In the left figure, 
    we sample $r_i$ from $\rmN(0,\sigma_i^2)$, where $\sigma_i^2$ is a sampled value using our Gibbs sampler, and 
    plot the $90\%$ credible interval for $|r_i|$.
    Similarly, in the right figure, we sample the error $\xi_i$ from $\rmN(0,\sigma_i^2 - \gamma^2)$ and plot the $90\%$ credible interval for $|\xi_i|$.
    For this test, we generated 29,500 posterior samples.
    }
    \label{fig:FN-R-sample-result}
\end{figure}
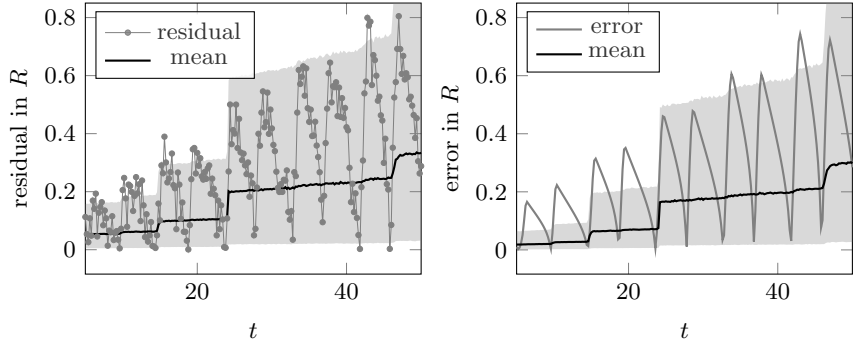

Figs.~\ref{fig:FN-V-sample-result} and~\ref{fig:FN-R-sample-result} plot a more direct quantification of the residual and error.
Instead of analyzing the posterior distribution of the discretization error variance $\sigma_i$, the left panels plot the corresponding posterior of the absolute residual $|r_i|$.
This is obtained by sampling $r_i$ from the normal distribution $\rmN(0,\sigma_i^2)$, where $\sigma_i^2$ is sampled using our Gibbs sampler.
The right panels, on the other hand, correspond to the errors, where $\rmN(0,\sigma_i^2-\gamma^2)$ is used instead of $\rmN(0,\sigma_i^2)$. 
As evident from the figures, a large proportion of data points fall within the 90$\%$ credible intervals.

\subsection{The Kepler equation}
We consider the Kepler equation
\begin{equation}
    \frac{\rmd}{\rmd t}
    \begin{bmatrix}
        q_1 \\ q_2 \\ p_1 \\ p_2
    \end{bmatrix}
    =
    \begin{bmatrix}
        p_1 \\ p_2 \\
        -\dfrac{q_1}{(q_1^2+q_2^2)^{3/2}}\\
        -\dfrac{q_2}{(q_1^2+q_2^2)^{3/2}}
    \end{bmatrix},
    \quad
    \begin{bmatrix}
        q_1(0) \\ q_2(0) \\ p_1(0) \\ p_2(0)
    \end{bmatrix}
    =
    \begin{bmatrix}
        1-e \\ 0 \\ 0 \\ \sqrt{(1+e)/(1-e)} 
    \end{bmatrix}
    =
    \begin{bmatrix}
        0.4 \\ 0 \\ 0 \\ 2
    \end{bmatrix}
    \quad (e=0.6),
\end{equation}
where $e$ denotes the eccentricity of the elliptic orbit.
For this example, we assume that we observe the velocity, i.e. $\calO(q_1,q_2,p_1,p_2) = \sqrt{p_1^2 + p_2^2}$ with the observation noise variance $\gamma^2 = 0.05^2$.
The velocity is observed from $t=10$ to $t=40$ with the interval $\Delta t = 0.2$
As a numerical integrator, we employ the symplectic Euler method with the step size $h=0.02$.

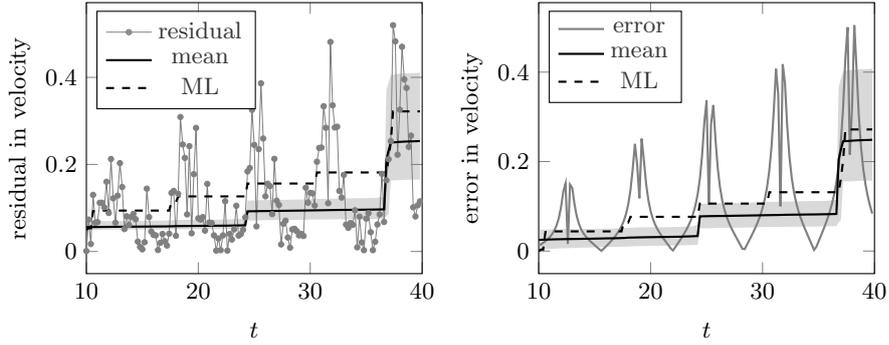
\begin{figure}
    \centering
    \begin{tikzpicture}
    \begin{axis}[
        width=6cm, 
        xlabel={$t$},
        ylabel={residual in velocity},
        ylabel style={yshift=-1em},
        legend pos=north west, 
        xmin=10, xmax=40,
        legend style={fill opacity=0.8, font=\small,align=left},
        tick label style={font=\small},
        label style={font=\small},
    ]
        \addplot[mark=*,  mark size=1pt, line width=0.1pt, gray] table [col sep=space] {Kepler/res_v.dat};
        \addlegendentry{residual}

        \addplot[mark=none, thick, black] table [col sep=space] {Kepler/mean_v.dat};
        \addlegendentry{mean}

        \addplot[mark=none, thick, dashed] table [col sep=space] {Kepler/iso_v.dat};
        \addlegendentry{ML}

        \addplot[name path =upper, draw=none] table [col sep=space] {Kepler/upper_v.dat};
        \addplot[name path =lower,draw=none] table [col sep=space] {Kepler/lower_v.dat};
        \addplot[fill opacity=0.3, fill=gray] 
        fill between[of=upper and lower];
    \end{axis}
\end{tikzpicture}
\begin{tikzpicture}

    \begin{axis}[
        width=6cm, 
        xlabel={$t$},
        ylabel={error in velocity},
        ylabel style={yshift=-1em},
        legend pos=north west, 
        xmin=10, xmax=40,
        legend style={fill opacity=0.8, font=\small,align=left},
        tick label style={font=\small},
        label style={font=\small},
    ]
        \addplot[thick, gray] table [col sep=space] {Kepler/err_v.dat};
        \addlegendentry{error}

        \addplot[mark=none, thick, black] table [col sep=space] {Kepler/err_mean_v.dat};
        \addlegendentry{mean}

        \addplot[mark=none, thick, dashed] table [col sep=space] {Kepler/err_iso_v.dat};
        \addlegendentry{ML}

        \addplot[name path =upper, draw=none] table [col sep=space] {Kepler/err_upper_v.dat};
        \addplot[name path =lower,draw=none] table [col sep=space] {Kepler/err_lower_v.dat};
        \addplot[fill opacity=0.3, fill=gray] 
        fill between[of=upper and lower];
    \end{axis}
\end{tikzpicture}

    \caption{Quantification results for the Kepler equation.
    The left figure shows the residual $r_i$, the quantification results in terms of $\sigma_i$ with their mean and $95\%$ credible interval, along with the quantification based on the maximum likelihood approach~\cite{mm21}.
    The right figure shows essentially the same information but with the observation removed: the error is plotted instead of the residual, and $\sqrt{\sigma_i^2 - \gamma^2}$ is plotted for the quantification results.}
    \label{fig:Kepler-result}
\end{figure}

\begin{figure}
    \centering
    \begin{tikzpicture}
    \begin{axis}[
        width=6cm, 
        xlabel={$t$},
        ylabel={residual in velocity},
        ylabel style={yshift=-1em},
        legend pos=north west, 
        xmin=10, xmax=40,
        legend style={fill opacity=0.8, font=\small,align=left},
        tick label style={font=\small},
        label style={font=\small},
    ]
        \addplot[mark=*,  mark size=1pt, line width=0.1pt, gray] table [col sep=space] {Kepler/res_v.dat};
        \addlegendentry{residual}

        \addplot[mark=none, thick, black] table [col sep=space] {Kepler/res_sample_mean.dat};
        \addlegendentry{mean}

        \addplot[name path =upper, draw=none] table [col sep=space] {Kepler/res_sample_upper.dat};
        \addplot[name path =lower,draw=none] table [col sep=space] {Kepler/res_sample_lower.dat};
        \addplot[fill opacity=0.3, fill=gray] 
        fill between[of=upper and lower];
    \end{axis}
\end{tikzpicture}
\begin{tikzpicture}

    \begin{axis}[
        width=6cm, 
        xlabel={$t$},
        ylabel={error in velocity},
        ylabel style={yshift=-1em},
        legend pos=north west,
        xmin=10, xmax=40,
        legend style={fill opacity=0.8, font=\small,align=left},
        tick label style={font=\small},
        label style={font=\small},
    ]
        \addplot[thick, gray] table [col sep=space] {Kepler/err_v.dat};
        \addlegendentry{error}

        \addplot[mark=none, thick, black] table [col sep=space] {Kepler/err_sample_mean.dat};
        \addlegendentry{mean}

        \addplot[name path =upper, draw=none] table [col sep=space] {Kepler/err_sample_upper.dat};
        \addplot[name path =lower,draw=none] table [col sep=space] {Kepler/err_sample_lower.dat};
        \addplot[fill opacity=0.3, fill=gray] 
        fill between[of=upper and lower];
    \end{axis}
\end{tikzpicture}

    \caption{
    These figures also present the residual and error plots for the velocity. 
    In the left figure, 
    we sample $r_i$ from $\rmN(0,\sigma_i^2)$, where $\sigma_i^2$ is a sampled value using our Gibbs sampler, and 
    plot the $90\%$ credible interval for $|r_i|$.
    Similarly, in the right figure, we sample the error $\xi_i$ from $\rmN(0,\sigma_i^2 - \gamma^2)$ and plot the $90\%$ credible interval for $|\xi_i|$.
    For this test, we generated 9,500 posterior samples.
    }
    \label{fig:Kepler-sample-result}
\end{figure}
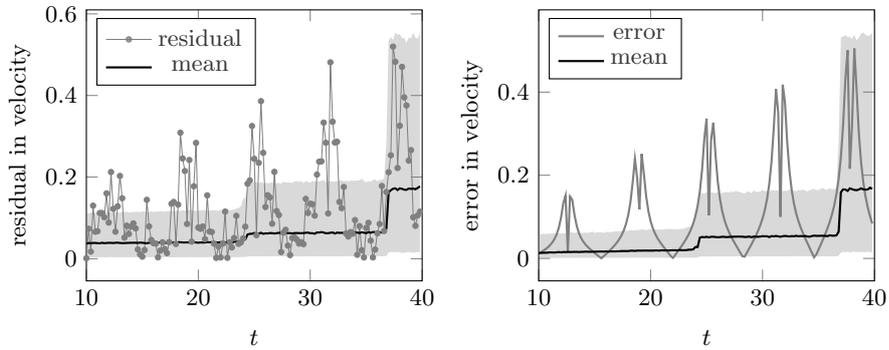

This is a test for the case where the observation operator is nonlinear.
The results are shown in Figs.~\ref{fig:Kepler-result} and~\ref{fig:Kepler-sample-result}, and trends similar to the FN case are observed.
The ratios of the residuals or errors included in the credible interval are different; however, the results capture the scale of propagation of the residual and error.

\section{Conclusion}\label{sec4}

In this paper, we have developed a Bayesian framework to quantify discretization errors in the numerical solutions of ordinary differential equations. 
The approach employs a shrinkage prior for the discretization error variances, and the Gibbs sampler can be used to sample from the corresponding posterior by approximating the $\log$-$\chi_1^2$ distribution with a Gaussian mixture model. 
Nevertheless, a significant challenge remains. 
It is anticipated that this methodology will eventually be incorporated into the parameter estimation process for ordinary differential equations.

\backmatter

\bibliography{references}

\end{document}